\documentclass[a4paper,11pt,twoside]{amsart}

\usepackage{amsmath}
\usepackage{amsthm}
\usepackage{amssymb}
\usepackage[all]{xypic}

\newtheorem{thm}{Theorem}[section]
\newtheorem{prop}[thm]{Proposition}
\newtheorem{lem}[thm]{Lemma}
\newtheorem{cor}[thm]{Corollary}

\theoremstyle{definition}
\newtheorem{definition}[thm]{Definition}
\newtheorem{remark}[thm]{Remark}
\newtheorem{ex}[thm]{Example}

\numberwithin{equation}{section}

\newcommand{\iso}{\cong}

\newcommand{\im}{\operatorname{im}}

\newcommand{\NS}{{\rm NS}}
\newcommand{\Pic}{{\rm Pic}}

\newcommand{\rk}{{\rm rk}\,}

\newcommand{\doublehat}[1]{\widehat{\widehat{\makebox[\width][l]{$#1$}}}  }

\newcommand{\cal}{\mathcal}

\newcommand{\kc}{{\cal C}}
\newcommand{\kd}{{\cal D}}

\newcommand{\kh}{{\cal H}}
\newcommand{\kk}{{\cal K}}

\newcommand{\kp}{{\cal P}}

\newcommand{\Db}{{\rm D}^{\rm b}}

\newcommand{\ZZ}{\mathbb{Z}}
\newcommand{\QQ}{\mathbb{Q}}
\newcommand{\RR}{\mathbb{R}}
\newcommand{\CC}{\mathbb{C}}

\newcommand{\PP}{\mathbb{P}}

\newcommand{\OO}{{\rm O}}

\begin{document}

\title[A finite group acting on the moduli space of K3
surfaces]{A finite group acting on the moduli space of K3 surfaces}

\author{Paolo Stellari}

\address{Dipartimento di Matematica ``F.
Enriques'', Universit{\`a} degli Studi di Milano, Via Cesare
Saldini 50, 20133 Milano, Italy}

\email{stellari@mat.unimi.it}

\keywords{K3 surfaces, moduli space of K3 surfaces, Fourier-Mukai
partners}

\subjclass[2000]{14J28, 14J10}

\begin{abstract} We consider the natural action of a finite group
on the moduli space of polarized K3 surfaces which induces a duality
defined by Mukai for surfaces of this type. We show that the group
permutes polarized Fourier-Mukai partners of polarized K3 surfaces
and we study the divisors in the fixed loci of the elements of this
finite group.\end{abstract}

\maketitle

\section{Introduction}\label{sec:intr}

In \cite{Mu1} Mukai introduced a duality for polarized K3
surfaces with Picard number 1. More precisely, given a K3 surface $X$ of this type, Mukai constructs
a fine moduli space $Y$ of stable sheaves on $X$ whose Chern classes
satisfy some restrictive conditions. Due to \cite{Mu2}, $Y$ is
also a K3 surface and, by \cite{Mu1}, if $X$ has a polarization of
degree $2d$ then $Y$ has a polarization of the same degree. In
analogy with the case of abelian varieties, $Y$ is interpreted as
a dual variety of $X$. The universal family on $X\times Y$ yields
an equivalence of the bounded derived categories of coherent
sheaves on $X$ and $Y$. Using the results in \cite{Mu2} it is quite easy to show that the same duality can be defined for K3 surfaces with arbitrary Picard group (Section \ref{subsec:Mudu}).

The aim of this paper is to show that Mukai's duality is given by
the natural action of a finite group $G_{2d}$ on the moduli space
$\kk_{2d}$ of K3 surfaces with a polarization of degree $2d$ ($d$ is
any positive integer). Indeed we prove that $G_{2d}$ acts by
interchanging periods of polarized K3 surfaces with equivalent
bounded derived categories of coherent sheaves. Moreover such an
action is free and transitive on the set of points of $\kk_{2d}$
parametrizing polarized K3 surfaces with the same derived category
and with Picard number 1. As we will point out, this immediately
leads to an interpretation of Mukai's duality as the action of
$G_{2d}$ on $\kk_{2d}$ (see Theorem \ref{thm:m1}).

The action of $G_{2d}$ is not necessarily free on the set of
periods of K3 surfaces with Picard number greater than 1. In
general, it seems difficult to determine explicitly the
irreducible components of the fixed loci and their dimension. Our
results on the geometrical properties of the divisors in the fixed
loci of the non-trivial elements of $G_{2d}$ are summarized by
Theorem \ref{thm:main1}. It essentially says that the irreducible
divisors contained in such fixed loci are Heegner divisors. As a
first approximation, a Heegner divisor can be thought of as the
image in $\kk_{2d}$ of a hyperplane section of the period domain
of polarized K3 surfaces and the K3 surfaces corresponding to its
points have Picard number greater or equal to $2$.

These divisors have recently appeared in other similar contexts. It is the case of \cite{Ni1} where the author studies rank-$2$ lattices $L$ such that all polarized K3 surfaces $X$ whose Picard group contains $L$ as a primitive sublattice are isomorphic to the moduli space of stable sheaves on $X$ determined by the choice of specific Chern classes (see Remark \ref{rmk:Nik} for more details).  

The quotient $G_{2d}\backslash\kk_{2d}$ parametrizes primitive
Hodge structures of K3 surfaces with polarizations of degree $2d$.
Such a variety naturally appears when describing divisors of the
moduli space of cubic 4-folds. A dense subset of
$G_{2d}\backslash\kk_{2d}$ parametrizes all triangulated
categories realized as derived categories of coherent sheaves on
K3 surfaces with Picard number 1 (see Section \ref{subsec:secint}
for these two remarks). As a consequence of Theorems \ref{thm:m1}
and \ref{thm:main1}, we describe the ramification divisor of the
natural map $\kk_{2d}\to G_{2d}\backslash\kk_{2d}$ and we show
that this dense subset does not meet the singular locus of
$G_{2d}\backslash\kk_{2d}$ (see Corollary \ref{cor:main1}).

\section{Mukai's duality and the action of $G_{2d}$}\label{sec:FMP}

Recall that a \emph{lattice} is a free
abelian group $L$ of finite rank with a non-degenerate symmetric
bilinear form $b:L \times L \rightarrow \mathbb{Z}$. We denote by
$\OO(L)$ the group of isometries of $L$. Given the \emph{dual
lattice} $L^\vee:=\mathrm{Hom}(L,\mathbb{Z})$ of a lattice $L$ and
the natural inclusion $L\hookrightarrow L^\vee$, the {\it
discriminant group of $L$} is the quotient $A_L:=L^\vee/L$. Notice
that $A_L$ inherits from $L$ a bilinear form with values in
$\QQ/\ZZ$. (For more details about lattices and discriminant groups,
see \cite{Ni}.) If $X$ is a K3 surface, the group
$H^2(X,\mathbb{Z})$ with the cup product is an even unimodular
lattice which is isomorphic to the lattice $\Lambda:=U^{\oplus
3}\oplus E_8(-1)^{\oplus 2}$ (for the meaning of $U$ and $E_8$ see
\cite{BPV}). Given two K3 surfaces $X$ and $Y$, a \emph{Hodge
isometry} is an isometry $\psi:H^2(X,\mathbb{Z})\rightarrow
H^2(Y,\mathbb{Z})$ such that $\psi(H^{2,0}(X))=H^{2,0}(Y)$. We
will denote by $\NS(X)$, $T(X):=\NS(X)^\perp$ and
$\rho(X):=\rk\NS(X)$ respectively the N{\'e}ron-Severi group,
the transcendental lattice and the Picard number of $X$.

The plan of this section is as follows. Mukai's duality is introduced in Section \ref{subsec:Mudu}. The
moduli space $\kk_{2d}$ of $2d$-polarized K3 surfaces and the action
of a finite group $G_{2d}$ on $\kk_{2d}$ are described in Section
\ref{sec:prel}. In that section we also discuss the interpretation
of Mukai's duality in terms of the action of $G_{2d}$. In Section
\ref{subsec:secint} we provide a few geometrical interpretations of
the quotient $G_{2d}\backslash\kk_{2d}$ which will be reconsidered
in Section \ref{sec:fixloc}.

\subsection{Mukai's duality}\label{subsec:Mudu} There exists a well-known duality for abelian varieties.
Indeed, given an abelian variety $A$ one can just consider the
abelian variety $\widehat A:=\Pic^0(A)$. The key aspects of this
duality are the following (see \cite{Mu}):
\begin{itemize}
\item[(A.1)] it is an involution, i.e.\
$\doublehat{A\:}=A$;\item[(A.2)] there exists an equivalence
$\Db(A)\iso\Db(\widehat A)$ induced by the Poincar\'e line bundle on
$A\times\widehat A$ (for a smooth projective variety $Y$, $\Db(Y)$
is the bounded derived category of coherent sheaves on
$Y$);\item[(A.3)] $A$ and $\widehat A$ have polarizations of the
same degree.
\end{itemize}

In \cite{Mu1}, Mukai proposed an analogous construction for K3
surfaces (i.e.\ 2-dimensional complex projective smooth varieties
with trivial canonical bundle and first Betti number $b_1=0$).
Recall that, for a positive integer $d$, a \emph{polarized K3
surface of degree $2d$} (or a \emph{$2d$-polarized K3 surface}) is
a pair $(X,\ell)$, where $X$ is a K3 surface and $\ell$ is a
primitive vector in the cone $C(X)^+\subset\Pic(X)\otimes\RR$
spanned by pseudo-ample (i.e.\ nef and big) divisors of $X$ and
$\ell^2=2d$.

Mukai's construction (\cite{Mu2,Mu1}) is as follows. Suppose that $(X,\ell)$ is a
$2d$-polarized K3 surface with $\ell^2=2d=2rs$ and $\gcd(r,s)=1$.
Consider the fine moduli space $M(r,\ell,s)$ parametrizing $\ell$-stable
sheaves $E$ on $X$ such that $\rk E=r$, $\mathrm{c}_1(E)=\ell$ and
$\chi(E)=r+s$. Deep results in \cite{Mu2} ensure that
$M(r,\ell,s)$ is a K3 surface as well. A K3 surface \emph{dual} to $X$ is
$\widehat X:=M(r,\ell,s)$. In complete analogy with the case of
abelian varieties, the results in \cite{Mu2,Mu1} imply that
\begin{itemize}
\item[(B.1)] Mukai's duality is an involution, i.e.\
$\doublehat{X}=X$ (see Remark \ref{rem:many}(iii) below);\item[(B.2)] there exists an equivalence
$\Db(X)\iso\Db(\widehat X)$ induced by the universal family on
$X\times\widehat X$;\item[(B.3)] $\widehat X$ has a polarization
$\ell'$ such that $(\ell')^2=2d$.
\end{itemize}

\begin{remark}\label{rem:many} (i) Mukai defined this duality in \cite{Mu1} just for K3 surfaces $X$ with Picard number 1. On the other hand it can be easily extended to polarized K3 surfaces $(X,\ell)$ with no resrictions on $\rho(X)$ using the explicit description in \cite{Mu2} of the lattice and Hodge structures of the second cohomology group of the moduli space $M(r,\ell,s)$, where $r$, $\ell$ and $s$ satisfy the previous conditions. Indeed, the only delicate part is (B.3). But due to \cite{Mu2}, if $v:=(r,\ell,s)$, then there exists a Hodge isometry $\psi:H^2(M(r,\ell,s),\ZZ)\to v^\perp/\ZZ v$ and $\psi^{-1}((0,\ell,2s))$ is a primitive vector in $\NS(M(r,\ell,s))$ with self-intersection $\ell^2$. 

(ii) Mukai's construction depends very much on the chosen factorization $2d=2rs$.
Indeed, \cite[Thm.\ 2.4]{St} shows that there may be many
non-isomorphic K3 surfaces dual (in the sense of Mukai) to a K3
surface.

(iii) Condition (B.1) can be better explained as follows: There exists a vector $\widehat{v}:=(\widehat{r},\widehat{\ell},\widehat{s})\in H^*(\widehat{X},\ZZ)$ with $\gcd(\widehat{r},\widehat{s})=1$, $\widehat{\ell}$ ample and primitive, $\widehat{\ell}^2=\ell^2=2\widehat{r}\widehat{s}$ and such that the fine moduli space $M(\widehat{r},\widehat{\ell},\widehat{s})$ of $\widehat{\ell}$-stable sheaves on $\widehat{X}$ (with the specified Chern classes) is isomorphic to $X$.\end{remark}

\subsection{The moduli space of $2d$-polarized K3
surfaces and the action of $G_{2d}$}\label{sec:prel} Let $h$ be a
primitive vector of $\Lambda$ with $h^2=2d$. The orthogonal
complement of $h$ in $\Lambda$ is isometric to the lattice
$$L_{2d}:=\langle k\rangle\oplus U^{\oplus 2}\oplus
E_8(-1)^{\oplus 2},$$ where $k^2=-2d$. The set
$D_{2d}:=\{\sigma\in\mathbb{P}(L_{2d}\otimes\mathbb{C}):
\sigma\cdot\sigma=0\mbox{ and }\sigma\cdot\overline{\sigma}>0\}$
is an open subset of a quadric in
$\mathbb{P}(L_{2d}\otimes\mathbb{C})$. We put $\Gamma(\langle
h\rangle):=\{g\in \OO(\Lambda):g(h)=h\}$ and
$\Gamma_{2d}:=\im(\Gamma(\langle h\rangle)\rightarrow
\OO(L_{2d}))$. It is easy to verify that $\Gamma_{2d}$ is the
kernel of the map $\Psi:\OO(L_{2d})\rightarrow\OO(A_{L_{2d}})$.
The quotient
\[
\kk_{2d}:=D_{2d}/\Gamma_{2d}
\]
is an irreducible quasi-projective variety of dimension $19$ (see
\cite{Ge}, Expos{\'e} XIII). By the Torelli Theorem and the
surjectivity of the period map, $\kk_{2d}$ is isomorphic to the
moduli space of $2d$-polarized K3 surfaces (see \cite{Ge,Do}).
Notice that by \cite[Thm.\ 1.14.4]{Ni}, the definition of $\kk_{2d}$
is independent of the choice of the primitive vector $h$ with self-intersection $2d$.

According to Section \ref{subsec:Mudu}, chosen two integers $r$ and
$s$ such that
\begin{eqnarray}\label{eqn:fond}
r,s\mbox{ positive}\;\;\;\;\;\;\;2d=2rs\;\;\;\;\;\;\gcd(r,s)=1,
\end{eqnarray}
Mukai's duality is given by the involution
\[
\delta_{r,s}:\kk_{2d}\longrightarrow\kk_{2d}
\]
sending the period $[\sigma_X]$ parametrizing a $2d$-polarized K3
surface $(X,\ell)$ to the period
$[\sigma_{r,s}]:=\delta_{r,s}([\sigma_X])$ parametrizing the
$2d$-polarized moduli space $(M(r,\ell,s),\ell')$. Therefore
Mukai's duality is realized by the action on $\kk_{2d}$ of the
group of automorphisms
\begin{eqnarray}\label{eqn:G2d}
\Delta_{2d}:=\langle\delta_{r,s}:r,s\mbox{ satisfying
\eqref{eqn:fond}}\rangle\subseteq\mathrm{Aut}(\kk_{2d})
\end{eqnarray}
generated by Mukai's involutions.

Consider the finite group
\[
G_{2d}:=\OO(A_{L_{2d}})/\{\pm \mathrm{id}\}.
\]
The main result of this section is the following theorem which will
be proved in Section \ref{subsec:proof} and which shows that Mukai's
duality can be described in terms of the natural action of $G_{2d}$.

\begin{thm}\label{thm:m1} The group $G_{2d}$ acts on $\kk_{2d}$ by sending the period of a
$2d$-polarized K3 surface $(X,\ell)$ to a period in the finite set
of $2d$-polarized Fourier-Mukai partners of $(X,\ell)$. Moreover
this action is free and transitive on the set of periods
corresponding to the Fourier-Mukai partners of any $2d$-polarized
K3 surface with Picard number 1 and $\Delta_{2d}=G_{2d}$.\end{thm}

Recall that two K3 surfaces $X_1$ and $X_2$ are
\emph{Fourier-Mukai partners} if there exists an exact equivalence
$\Db(X_1)\cong \Db(X_2)$. We write $\mathrm{FM}(X)$
for the set of isomorphism classes of Fourier-Mukai partners of
$X$.

Before passing to the proof of Theorem \ref{thm:m1}, we fix some
notations which will be used for the rest of this paper. By
definition, there exists a (surjective) map $\OO(L_{2d})\to G_{2d}$.
Hence given $f\in\OO(L_{2d})$, we write $[f]$ for the image
of $f$ in $G_{2d}$. If $\sigma$ is an element of $D_{2d}$ and
$\pi:D_{2d}\rightarrow\kk_{2d}$ is the natural projection,
$[\sigma]:=\pi(\sigma)$. If $\sigma\in D_{2d}$, then $T(\sigma)$ is
the minimal primitive sublattice of $\Lambda$ such that $\sigma\in
T(\sigma)\otimes\mathbb{C}$. We put
\[
\rho([\sigma]):=22-\rk T(\sigma)\;\;\;\;\;\mbox{and}
\;\;\;\;\;\OO_\sigma:=\{g\in
\OO(T(\sigma)):g(\mathbb{C}\sigma)=\mathbb{C}\sigma\}.
\]
Observe that if $[\sigma]=[\tilde{\sigma}]$ then $\OO_\sigma$ is
conjugate to $\OO_{\tilde{\sigma}}$. Moreover, we define
$(X_\sigma,\ell_\sigma)$ to be the $2d$-polarized K3 surface
corresponding to $[\sigma]$ (in particular $T(X_\sigma)\iso
T(\sigma)$). The points of $\kk_{2d}$ will be often called
periods.

\subsection{Proof of Theorem \ref{thm:m1}}\label{subsec:proof} Let
us start from the following easy lemma.

\begin{lem}\label{lem:act} The group $G_{2d}$ acts on $\kk_{2d}$ and
$G_{2d}\cong(\mathbb{Z}/2\mathbb{Z})^{p(d)-1}$ where $p(1)=1$ and
$p(d)$ is the number of distinct primes dividing $d$, if $d\geq
2$.\end{lem}

\begin{proof} Obviously $\OO(L_{2d})$ acts on $D_{2d}$. The subgroup $\Gamma_{2d}$ is the kernel of the map
$\Psi:\OO(L_{2d})\rightarrow \OO(A_{L_{2d}})$. Hence it is normal in
$\OO(L_{2d})$ and $\OO(L_{2d})/\Gamma_{2d}$ acts on
$\kk_{2d}=\Gamma_{2d}\backslash D_{2d}$. Since, under our
hypothesis, $\Psi$ is onto (see \cite[Thm.\ 1.14.2]{Ni}),
$\OO(A_{L_{2d}})\cong \OO(L_{2d})/\Gamma_{2d}$ acts on $\kk_{2d}$.
As $-\mathrm{id}$ acts trivially on $D_{2d}$, $G_{2d}$ acts on
$\kk_{2d}$. The second part of the statement is \cite[Lemma
3.6.1]{Sc} (see also \cite{Og}).\end{proof}

The fact that $G_{2d}$ acts on $\kk_{2d}$ by sending the period of
a $2d$-polarized K3 surface $(X,\ell)$ to the period of a
$2d$-polarized Fourier-Mukai partner of $(X,\ell)$ is easy.
Indeed, given $[\sigma]\in \kk_{2d}$ corresponding to $(X,\ell)$
and an $f\in \OO(L_{2d})$, then $f(T(\sigma))=T(f(\sigma))$. If
$(X',\ell')$ corresponds to $[f(\sigma)]\in\kk_{2d}$, then there
is a Hodge isometry between $T(X)$ and $T(X')$. To conclude,
observe that Orlov proved in \cite{Or} (using results of Mukai)
that, given two K3 surfaces $X$ and $Y$, the following three
conditions are equivalent: (i)
$\mathrm{D^b}(X)\cong\mathrm{D^b}(Y)$; (ii) there exists a Hodge
isometry $T(X)\cong T(Y)$; (iii) $Y$ is isomorphic to a smooth
compact 2-dimensional fine moduli space of stable sheaves on $X$.

\begin{definition}\label{def:FM} Let $(X,\ell)$ be a
$2d$-polarized K3 surface. We define the set of isomorphism classes
of $2d$-polarized Fourier-Mukai partners of $(X,\ell)$:
\[
\mathrm{FM}_{2d}(X,\ell):=\left\{(X',\ell'):\begin{array}{l}
\;(1)\;\; (X',\ell')\mbox{ is a $2d$-polarized  K3 surface} \\
\;(2)\;\;
\mathrm{D^b}(X)\cong\mathrm{D^b}(X')\end{array}\right\}/\cong,
\]
where $(X_1,\ell_1)\iso(X_2,\ell_2)$ if and only if there exists
an isomorphism $\psi:X_1\rightarrow X_2$ such that
$\psi^*(\ell_2)=\ell_1$.\end{definition}

\begin{lem}\label{prop:finite} Let $(X,\ell)$ be a $2d$-polarized K3 surface. Then the set
$\mathrm{FM}_{2d}(X,\ell)$ is finite.\end{lem}

\begin{proof} By \cite[Prop.\ 5.3]{BM}, the set $\mathrm{FM}(X)$ is
finite. Let $m:=|\mathrm{FM}(X)|$ and let $Y_1,\cdots, Y_m$ be
representatives of the isomorphism classes of the (unpolarized)
Fourier-Mukai partners of $X$.

We first show that the number of non-isomorphic $2d$-polarizations
on $Y_j$, for $j\in\{1,\ldots,m\}$, is finite. Indeed, consider the
set $A:=\{c\in C(Y_j)^+:c^2=2d\}$. Given $c\in A$, by Saint-Donat's
result in \cite{Sa} and Bertini's Theorem, there exists an
irreducible $D$ in the linear system $|3c|$ such that $D^2=18d$. If
$$B:=\{|D|:D\in\mathrm{Div}(Y_j)\mbox{ is irreducible and
}D^2=18d\},$$ item (b) of Theorem 0.1 in \cite{Sr}, asserts that
$B/\mathrm{Aut}(Y_j)$ is finite. Since there exists an injective
map $A\hookrightarrow B$ defined by $c\mapsto |3c|$, $Y_j$ has
only a finite number $P_d(Y_j)$ of non-isomorphic
$2d$-polarizations. It is easy to see that
$$|\mathrm{FM}_{2d}(X,\ell)|=\sum_{j=1}^m P_d(Y_j).$$ Hence
$\mathrm{FM}_{2d}(X,\ell)$ is finite.\end{proof}

Consider the dense subset $$\kc_\mathrm{gen}\subseteq\kk_{2d}$$
parametrizing $2d$-polarized K3 surfaces with Picard number 1.

\begin{lem}\label{prop:acFM} If $(X,\ell)$ is a $2d$-polarized K3 surface
corresponding to $[\sigma]\in\kc_\mathrm{gen}$, then the periods of
distinct isomorphism classes in $\mathrm{FM}_{2d}(X,\ell)$
correspond to distinct points in the orbit $G_{2d}\cdot[\sigma]$
and $|G_{2d}\cdot[\sigma]|=|\mathrm{FM}_{2d}(X,\ell)|$.\end{lem}

\begin{proof} Take $(X,\ell)$ and $[\sigma]\in\kk_{2d}$ as in the statement.
In particular, $\mathrm{NS}(X)\iso\langle h\rangle$ (here $h$ is a
generator of the orthogonal complement of $L_{2d}$ in $\Lambda$). By
\cite[Thm.\ 1.14.4]{Ni}, $\langle h\rangle$ and $\mathrm{NS}(X)$
have unique (up to isometries of $\Lambda$) primitive embeddings in $\Lambda$. Therefore,
$$\mathrm{FM}_{2d}(X,\ell)=\{X_{f(\sigma)}:f\in
\OO(L_{2d})\}/\mathrm{isom}.$$ Given $f\in \OO(L_{2d})$, by the
Torelli Theorem, $X_{f(\sigma)}\cong X_\sigma$ if and only if either
$f\in \OO_\sigma$ or $f$ extends to a Hodge isometry in
$\OO(\Lambda)$. Let
$$I_1:=\im(\OO(\langle h\rangle)\rightarrow \OO(A_{\langle h\rangle})\cong
\OO(A_{L_{2d}}))\;\;\;\;\;\mbox{and}\;\;\;\;\;I_2:=\im(\Psi:\OO(L_{2d})\rightarrow
\OO(A_{L_{2d}})).$$Obviously, $I_1=\{\pm \mathrm{id}\}$ and, by
\cite[Thm.\ 1.14.2]{Ni}, $I_2=\OO(A_{L_{2d}})$. Now \cite[Lemma
4.1]{Og} implies that the image of the natural composition
$\OO_\sigma\rightarrow \OO(L_{2d})\rightarrow \OO(A_{L_{2d}})$ is
$\{\pm \mathrm{id}\}$.

The results \cite[Prop.\ 1.6.1, Cor.\ 1.5.2]{Ni} state that $g\in
\OO(L_{2d})$ lifts to an isometry in $\OO(\Lambda)$ if and only if
$\Psi (g)\in \{\pm \mathrm{id}\}$. In particular, if $g,h\in G_{2d}$
and $g\neq h$, then $g([\sigma])\neq h([\sigma])$. Hence
$\{g([\sigma]):g\in G_{2d}\}$ is the set of periods of the
non-isomorphic Fourier-Mukai partners of $X$. Hence we have just
shown that distinct isomorphism classes in
$\mathrm{FM}_{2d}(X,\ell)$ correspond to distinct points in the
orbit $G_{2d}\cdot[\varphi(\sigma_X)]$ and that
$$|G_{2d}\cdot[\varphi(\sigma_X)]|=|\mathrm{FM}_{2d}(X,\ell)|.$$ (Compare
this with the description of $G_{2d}$ and of its order given in
Lemma \ref{lem:act}.)\end{proof}

This concludes the proof of the first part of Theorem
\ref{thm:m1}. Now we just need to show that $\Delta_{2d}=G_{2d}$.

Let $(X,\ell)$ be a $2d$-polarized K3 surface corresponding to
$[\sigma_{(X,\ell)}]\in\kc_\mathrm{gen}$. Due to \cite[Thm.\
2.4]{St}, for any Fourier-Mukai partner $Y$ of $X$, there are two
positive integers $r$ and $s$ satisfying \eqref{eqn:fond} and such
that $Y$ is isomorphic to the fine moduli space $M(r,\ell,s)$
parametrizing stable sheaves $E$ on $X$ with $\rk E=r$, ${\rm
c}_1(E)=\ell$ and $\chi(E)-\rk E=s$. Moreover, $M(r,\ell,s)$ has
a (unique) $2d$-polarization $\ell'_X$. We denote by
$[\sigma_{(M(r,\ell,s),\ell'_X)}]\in\kc_\mathrm{gen}$ the period
corresponding to the $2d$-polarized K3 surface
$(M(r,\ell,s),\ell'_X)$.

Since $G_{2d}$ acts by sending the period of a $2d$-polarized K3
surface $(X,\ell)$ to the period of a $2d$-polarized Fourier-Mukai
partner of $(X,\ell)$, by \cite[Thm.\ 1.4]{Mu2}, given $g\in
G_{2d}$ there exist two integers $r$ and $s$ satisfying
\eqref{eqn:fond} and such that
\[
g:[\sigma_{(X,\ell)}]\longmapsto
[\sigma_{(M(r,\ell,s),\ell'_X)}],
\]
for any $[\sigma_{(X,\ell)}]\in\kc_\mathrm{gen}$. In other words,
\[
\delta_{r,s}|_{\kc_\mathrm{gen}}=g|_{\kc_\mathrm{gen}}.
\]
Conversely, given $r$ and $s$ as in \eqref{eqn:fond} and
$[\sigma_{(X,\ell)}]\in\kc_\mathrm{gen}$, by Lemma \ref{prop:acFM}
there is a $g\in G_{2d}$ such that
\[
g([\sigma_{(X,\ell)}])=[\sigma_{(M(r,\ell,s),\ell'_X)}].
\]
Reasoning as before, this yields the equality
$\delta_{r,s}|_{\kc_\mathrm{gen}}=g|_{\kc_\mathrm{gen}}$.

To conclude the proof of Theorem \ref{thm:m1}, observe that since
$\kc_\mathrm{gen}$ is dense in $\kk_{2d}$, $\delta_{r,s}=g$ whenever
$\delta_{r,s}|_{\kc_\mathrm{gen}}=g|_{\kc_\mathrm{gen}}$.

\subsection{A geometrical interpretation}\label{subsec:secint} Consider the quotient
\[
\kp_{2d}:=G_{2d}\backslash\kk_{2d}.
\]
It is very easy to see that $\kp_{2d}=\OO(L_{2d})\backslash
D_{2d}$. Hence, due to the surjectivity of the period map and the
Torelli Theorem for K3 surfaces, each point of $\kp_{2d}$
naturally parametrizes all K3 surfaces with Hodge-isometric primitive (with
respect to a polarization of degree $2d$) second cohomology groups.
Recall that for a $2d$-polarized K3 surface $(X,\ell)$, its
primitive second cohomology group is the
orthogonal complement of $\ell$ in $H^2(X,\ZZ)$. Such a sublattice
naturally inherits form $H^2(X,\ZZ)$ a weight-two Hodge structure.

Very often the quotient $\kp_{2d}$ can be seen as a birational model
of special divisors of the moduli space $\kc$ of cubic 4-folds (i.e.\ hypersurfaces of degree $3$ in $\PP^5$). This
is the case if and only if $2d$ is not divisible by $4$, $9$ or any
odd prime $p\equiv 2\pmod{3}$. Indeed, under these assumptions, one
considers, as in \cite{H2}, the divisors $\kc_{2d}$ of $\kc$ whose
points parametrize cubic 4-folds $X$ such that the lattice
$$A(X):=H^4(X,\ZZ)\cap H^{2,2}(X)$$
contains a primitive positive definite sublattice $K_{2d}$ with the
properties:
\begin{itemize}
\item the class $H^2$ (here $H$ is the hyperplane section of $X$)
belongs to $K_{2d}$; \item $\rk K_{2d}=2$ and $2d=|A_{K_{2d}}|$.
\end{itemize}
Interesting examples of divisors of this type are the ones whose
generic points parametrize cubic 4-folds $X$ whose Fano variety
$F(X)$ of lines contained in $X$ is birational to the Hilbert scheme
$S^{[2]}$, where $S$ is a K3 surface. Now \cite[Thm.\ 1.0.2]{H2}
allows us to conclude the existence of a birational map
$\kc_{2d}\stackrel{\sim}{\dashrightarrow}\kp_{2d}$.

Let us investigate a little bit more the geometry of $\kp_{2d}$. In
particular, consider the dense subset
$\kc_\mathrm{gen}\subset\kk_{2d}$ whose points are the periods of K3
surfaces with Picard number 1. Due to Theorem \ref{thm:m1}, one can
consider the dense subset
$$\widetilde\kc_\mathrm{gen}:=G_{2d}\backslash\kc_\mathrm{gen}\subset\kp_{2d}.$$
Applying once more Theorem \ref{thm:m1} we see that
$\widetilde\kc_\mathrm{gen}$ parametrizes all triangulated
categories realized as bounded derived categories of coherent
sheaves on K3 surfaces with Picard number 1 and of degree $2d$.

\section{Fixed loci of the action of $G_{2d}$ and divisors of
$\kk_{2d}$}\label{sec:fixloc}

By Theorem \ref{thm:m1}, Mukai's duality has no fixed points
contained in the dense subset $\kc_\mathrm{gen}$. On the other
hand it may happen that some $[\sigma]$ in the complement of
$\kc_\mathrm{gen}$ is fixed by a non-trivial
$\delta\in\Delta_{2d}$. In this section we study the divisors
contained in the fixed loci of non-trivial elements of
$\Delta_{2d}$ or, equivalently, of $G_{2d}$.

Proceeding in this direction we state a few definitions. Given a
primitive sublattice $L\hookrightarrow L_{2d}$ with signature
$(2,18)$, we put
\[
D_{L}:=\{\sigma\in\mathbb{P}(L\otimes\mathbb{C}):\sigma^2=0\mbox{
and
}\sigma\cdot\overline{\sigma}>0\}=\mathbb{P}(L\otimes\mathbb{C})\cap
D_{2d}.
\]
Such a primitive sublattice defines a divisor
$\kh_{L}:=\im(\pi:D_{L^\perp}\longrightarrow \kk_{2d})$, where
$\pi$ is the restriction to $D_{L^\perp}$ of the natural
projection $D_{2d}\rightarrow \kk_{2d}$. Obviously, if
$j\in\Gamma_{2d}$, then $\kh_L=\kh_{j(L)}$.

\begin{definition} A divisor $\kd$ of $\kk_{2d}$ is a \emph{Heegner divisor}
if there exists a primitive sublattice $L$ of $L_{2d}$ with
signature $(2,18)$ such that $\kd=\kh_{L}$.\end{definition}

\begin{remark}\label{rmk:Nik} Heegner divisors appear also in \cite{Ni1}. Indeed Nikulin studies primitive sublattices $L$ of $\Lambda$ such that, fixed two positive integers $r,s$ and $l\in L$ ($l^2>0$),
\begin{itemize}
\item[(1)] for all polarized K3 surfaces $(X,\ell)$ with a primitive embedding $i:L\hookrightarrow\NS(X)$ and $\ell=i(l)$, there is an isomorphism $X\iso M(r,\ell,s)$;

\item[(2)] property (1) does not hold true if we substitute $L$ with a primitive sublattice $L'\hookrightarrow L$ with $\rk L'<\rk L$ and $l\in L'$.
\end{itemize}
Following Nikulin's terminology, such a lattice $L$ is called \emph{critical polarized K3 Picard lattice}. It is clear that a critical polarized K3 Picard lattice of rank $2$ defines a Heegner divisor of the moduli space  $\kk_{l^2}$ of $l^2$-polarized K3 surfaces. Moreover, the generic points of these divisors (i.e.\ the periods in $\kk_{l^2}$ parametrizing K3 surfaces $X$ with $\NS(X)\iso L$) are fixed points of Mukai's duality if $\gcd(r,s)=1$ (simply use Theorem \ref{thm:m1}).\end{remark}

For a Heegner divisor $\kh_L$ we denote by $V(\kh_L)$ the orthogonal
complement of $L$ in $L_{2d}$. If $\beta\in L_{2d}$, we write
$$\beta=\alpha_\beta k+m_\beta j_\beta\in L_{2d}=\langle
k\rangle\oplus U^{\oplus 2}\oplus E_8(-1)^{\oplus 2},$$ where
$j_\beta\in U^{\oplus 2}\oplus E_8(-1)^{\oplus 2}$ is primitive and
$\alpha_\beta,m_\beta\in\mathbb{Z}$. Consider the set
\begin{eqnarray}\label{eqn:refl}
\Re:=\left\{\langle\beta\rangle\subset L_{2d}:\begin{array}{l}
(\mathrm{a})\;\;\mbox{$\beta\in L_{2d}$ is primitive and $\beta^2$
divides $2m_{\beta}$ and $2d$}
\\ (\mathrm{b})\;\;\mbox{$\beta^2<0$, $\beta^2\neq -2$ and $\beta^2\neq -2d$}\end{array}\right\}.
\end{eqnarray}

In Sections \ref{subsec:fixlo} and \ref{subsec:nonempty} we will
prove the following result:

\begin{thm}\label{thm:main1} The irreducible divisors of $\kk_{2d}$ contained in the
union of the fixed loci of the non-trivial elements of $G_{2d}$ are
the Heegner divisors $\kh_L$ such that $V(\kh_L)\in\Re$.\end{thm}

As a consequence of Theorems \ref{thm:m1} and \ref{thm:main1} we can
prove the following:

\begin{cor}\label{cor:main1} The irreducible components of the ramification divisor of the natural map
\[
p:\kk_{2d}\longrightarrow\kp_{2d}
\]
are the Heegner divisors $\kh_L$ such that $V(\kh_L)\in\Re$.
Moreover $\widetilde\kc_\mathrm{gen}$ is contained in the smooth
locus of $\kp_{2d}$.\end{cor}

\begin{proof} The first part is just Theorem \ref{thm:main1}.
For the second part observe that, due to Theorem \ref{thm:m1},
$\mathrm{Sing}(\kp_{2d})\cap\widetilde\kc_\mathrm{gen}=\mathrm{Sing}(\kk_{2d})\cap\kc_\mathrm{gen}$,
where $\mathrm{Sing}(\kp_{2d})$ and $\mathrm{Sing}(\kk_{2d})$ are
the singular loci of $\kp_{2d}$ and $\kk_{2d}$ respectively. But
given $[\sigma]\in\kk_{2d}$ corresponding to a $2d$-polarized K3
surface $(X_\sigma,\ell_\sigma)$ with $\rho(X_\sigma)=1$,
$\OO_\sigma=\{\pm \mathrm{id}\}$ (see \cite[Lemma 4.1]{Og}). Hence
$[\sigma]\not\in\mathrm{Sing}(\kk_{2d})$ and this concludes the
proof.\end{proof}

\subsection{Fixed loci and reflections}\label{subsec:fixlo} As a first step towards the proof
of Theorem \ref{thm:main1}, we show that a divisor in the fixed
locus of a non-trivial $g\in G_{2d}$ is union of Heegner divisors.
Moreover, if $g\in G_{2d}$ is non-trivial and fixes a divisor then
it is representable by a reflection.

\begin{definition}\label{def:dr} (i) If $g\in G_{2d}$, the set $\mathrm{Fix}(g):=\{[\sigma]\in
\kk_{2d}:g([\sigma])=[\sigma]\}$ is the \emph{fixed locus of $g$}.

(ii) Let $\beta\in L_{2d}$ be such that $\beta^2\neq 0$ and
$\frac{2x\cdot\beta}{\beta^2}\beta\in L_{2d}$, for every $x\in
L_{2d}$. Then the isometry $r_\beta\in\OO(L_{2d})$ defined by
$$r_\beta(x):=x-2\frac{x\cdot\beta}{\beta^2}\beta,$$ for any $x\in
L_{2d}$, is the \emph{reflection with respect to $\beta$}.

(iii) The subset of $G_{2d}$ containing the non-trivial elements
representable by reflections with respect to vectors with negative
self-intersection is denoted by
\[
\mathrm{Ref}(G_{2d}):=\{g\in
G_{2d}\setminus\{\mathrm{id}\}:g=[r_\beta], \mbox{ with
}\beta^2<0\mbox{ and }\beta\in L_{2d}\}.
\]\end{definition}

Observe that, due to Theorem \ref{thm:m1}, if $[\sigma]\in\kk_{2d}$
is fixed by some non-trivial $g\in G_{2d}$, then $\rho([\sigma])\geq
2$. Under this additional requirement, we can prove the following:

\begin{lem}\label{prop:m1} A point $[\sigma]\in\kk_{2d}$ with $\rho([\sigma])=2$ is fixed
by at most one  non-trivial $g\in G_{2d}$.\end{lem}

\begin{proof} If $|G_{2d}|=1$ the
result is trivial. So, we suppose $|G_{2d}|>1$.

It is very easy to observe that if $[\sigma]\in \kk_{2d}$, then
$[\sigma]\in \mathrm{Fix}(g)$ if and only if there exists an $f\in
\OO(L_{2d})$ such that $[f]=g$, $f(T(\sigma))\subseteq
T(\sigma)$ and $f|_{T(\sigma)}\in \OO_{\sigma}$. Thus the
proposition is proved if, given
\[
S_\sigma:=\left\{f\in \OO(L_{2d}): \begin{array}{l}
(\mathrm{a})\;\; f(T(\sigma))\subseteq T(\sigma) \\
(\mathrm{b})\;\; f|_{T(\sigma)}\in \OO_{\sigma}\end{array}\right\},
\]
we show that $|\Xi(S_\sigma)|\leq 2$, where $\Xi$ is the
composition of $\Psi:\OO(L_{2d})\rightarrow\OO(A_{L_{2d}})$ with
the natural projection from $\OO(A_{L_{2d}})$ onto $G_{2d}$.

As $\rho([\sigma])=2$, there is a primitive $\beta\in L_{2d}$
(unique up to sign) such that
$(L_{2d})_\mathbb{Q}=\langle\beta\rangle_\mathbb{Q}\oplus
T(\sigma)_\mathbb{Q}$, where, for a lattice $L$, we write
$L_\mathbb{Q}:=L\otimes_\mathbb{Z}\mathbb{Q}$. This decomposition
of $(L_{2d})_\mathbb{Q}$ induces an injective homomorphism
$$S_\sigma\hookrightarrow\OO(\langle\beta\rangle)\times\OO(T(\sigma)).$$
From Proposition B.1 in \cite{HLOY1} it follows that
$\im(S_\sigma\rightarrow\OO(T(\sigma)))$ is a cyclic group of order
$2n$. Moreover, $\OO(\langle\beta\rangle)=\{\pm \mathrm{id}\}$.
Hence $S_\sigma$ can be generated by at most two elements and so
$|\Psi(S_\sigma)|\leq 4$ (indeed, due to Lemma \ref{lem:act}, for any $g\in G_{2d}$,
$g^2=\mathrm{id}$). Since $-\mathrm{id}_{L_{2d}}\in\OO(L_{2d})$ lies
in $S_\sigma$ and it maps to $(-\mathrm{id},-\mathrm{id})$ in
$\OO(\langle\beta\rangle)\times\OO(T(\sigma))$, $\Xi(S_\sigma)$ can
be generated by one element and thus $|\Xi(S_\sigma)|\leq
2$.\end{proof}

\begin{remark}\label{lem:zero} Assume that $[\sigma]\in\kk_{2d}$ is a zero-dimensional irreducible component
of $\mathrm{Fix}(g)$, for some non-trivial $g\in G_{2d}$. Then
$\rho([\sigma])$ is even. Indeed, generalizing the proof of
\cite[Lemma 4.1]{Og}, we see that for each period $[\sigma]$
corresponding to a K3 surface with odd Picard number,
$\OO_\sigma=\{\pm \mathrm{id}\}$. Thus if $\rho([\sigma])$ is odd
and $[\sigma]\in \mathrm{Fix}(g)$, then
$\kh_{T(\sigma)}:=\im(D_{T(\sigma)}\rightarrow\kk_{2d})$ is
contained in $\mathrm{Fix}(g)$ and it is a subvariety of codimension
$0<\rho([\sigma])-1<19$ in $\kk_{2d}$. As
$[\sigma]\in\kh_{T(\sigma)}$, it follows that $[\sigma]$ is not a
zero-dimensional irreducible component of
$\mathrm{Fix}(g)$.\end{remark}

Now we consider the irreducible divisors of $\kk_{2d}$ contained
in the fixed loci.

\begin{prop}\label{prop:m2} A divisor
$\mathcal{D}$ of $\kk_{2d}$ is fixed by at most one non-trivial
$g\in G_{2d}$. In this case $g\in \mathrm{Ref}(G_{2d})$ and the
irreducible components of $\mathcal{D}$ are Heegner divisors.
Conversely, if $g\in \mathrm{Ref}(G_{2d})$, then $\mathrm{Fix}(g)$
contains a Heegner divisor.\end{prop}

\begin{proof} Obviously, $[\sigma]\in\mathcal{D}\subseteq\mathrm{Fix}(g)$
for a non-trivial $g\in G_{2d}$ if and only if there exists an $f\in
\OO(L_{2d})$ such that $[f]=g$ and $\sigma$ is an eigenvector of
the $\mathbb{C}$-linear extension $f_\mathbb{C}$ of $f$ for a
certain complex eigenvalue $\lambda$. More precisely, $\lambda=\pm
1$ otherwise the eigenspace of $\lambda$ would have dimension less
or equal to 10 (indeed, $\overline{\lambda}$ is also an eigenvalue of
$f_\CC$).

This and the fact that, due to Theorem \ref{thm:m1}, a point
$[\sigma]$ such that $\rho([\sigma])=1$ is not fixed, imply that
there exist a subset $S\subseteq \{f\in
\OO(L_{2d}):[f]=g\}$ and, for every $f\in S$, a non-trivial
$\beta_f\in L_{2d}$ such that either
$f_\CC|_{\beta_f^\perp}=\mathrm{id}$ or
$f_\CC|_{\beta_f^\perp}=-\mathrm{id}$ and
\begin{eqnarray}\label{eqn:2}
\mathcal{D}=\pi\left(\bigcup_{f\in
S}\mathbb{P}(\beta_f^\perp\otimes\mathbb{C})\cap D_{2d}\right),
\end{eqnarray}
where $\pi:D_{2d}\rightarrow\kk_{2d}$ is the usual projection. This
means that $\mathcal{D}$ is the union of the Heegner divisors
$\{\kh_{\beta_f^\perp}\}_{f\in S}$.

A Heegner divisor $\kh_{L}$ of $\kk_{2d}$, where $L$ is a primitive
sublattice of $L_{2d}$ whose signature is $(2,18)$, is irreducible.
Indeed, by \cite[Cor.\ 2.5]{Mo}, there is a primitive embedding
$M:=L^\perp\hookrightarrow U\oplus U$ and hence a primitive
embedding $U\hookrightarrow M^\perp$. Take $D_{M^\perp}$ and the
groups $\Gamma(M):=\{g\in \OO(\Lambda):g(m)=m\mbox{ for any }m\in
M\}$ and $\Gamma_M:=\im(\Gamma(M)\rightarrow \OO(M^\perp))$. By
\cite[Prop.\ 5.6]{Do}, the quotient $\kk_M:=D_{M^\perp}/\Gamma_M$ is
irreducible. Since $\Gamma(M)\subseteq\Gamma(\langle h\rangle)$,
there exists a surjective map $\kk_{M}\rightarrow\kh_L$ and thus
$\kh_L$ is irreducible. From this we conclude that the irreducible
components of $\kd$ are Heegner divisors.

A straightforward consequence of the previous remarks and of
(\ref{eqn:2}) is that, when $f\in S$, $f=r_{\beta_f}$ because
$f|_{\beta_f^\perp}$ is $\pm\mathrm{id}$. Thus if a Heegner divisor
$\kh_{L}$ is contained in $\mathrm{Fix}(g)$, for $g\in G_{2d}$
non-trivial, then $g\in\mathrm{Ref}(G_{2d})$. Now let us suppose
that $\kh_{L}\subset\mathrm{Fix}(g_1)\cap\mathrm{Fix}(g_2)$, for
$g_1,g_2\in G_{2d}$ non-trivial. Lemma \ref{prop:m1}, applied to
$[\sigma]\in\kh_{L}$ with $\rho([\sigma])=2$, shows that $g_1=g_2$.
Hence, any divisor $\mathcal{D}$ is fixed by at most one non-trivial
$g\in G_{2d}$.

The last statement of the proposition is trivial.\end{proof}

\subsection{Non-empty fixed loci}\label{subsec:nonempty} In this section we conclude the proof of Theorem \ref{thm:main1},
describing the Heegner divisors fixed by the non-trivial elements of
$G_{2d}$ .

\begin{lem}\label{thm:Heegnerfixed} If
$L$ is a primitive sublattice of $L_{2d}$ with signature $(2,18)$,
then the Heegner divisor $\kh_{L}$ is fixed by a non-trivial $g\in
G_{2d}$ if and only if $V(\kh_L)\in\Re$.\end{lem}

\begin{proof} A Heegner divisor
$\kh_{L}$ with $V(\kh_L)=\langle\beta\rangle$, for some primitive
$\beta\in\ L_{2d}$ is fixed by a non-trivial $g\in G_{2d}$ if and
only if $\beta$ defines a reflection $r_{\beta}$ such that
$[r_{\beta}]=g\in\mathrm{Ref}(G_{2d})$ (see Proposition
\ref{prop:m2}).

Suppose that $\beta$ defines a reflection in $\OO(L_{2d})$ which is
non-trivial in $G_{2d}$. As at the beginning of Section
\ref{sec:fixloc}, let us write $\beta=\alpha k+m j\in\langle
k\rangle\oplus U^{\oplus 2}\oplus E_8(-1)^{\oplus 2}$, where $j\in
U^{\oplus 2}\oplus E_8(-1)^{\oplus 2}$ is primitive and
$\alpha,m\in\mathbb{Z}$. Let $-\beta^2/2=p^eb_1$, with $(b_1,p)=1$ and $p$ prime.
Since $j\in U^{\oplus 2}\oplus E_8(-1)^{\oplus 2}$ and $j\neq 0$
(otherwise $[r_\beta]=\mathrm{id}\in G_{2d}$) there exists
an $x\in L_{2d}$ such that $x\cdot\beta=m$. Since
$\frac{2(x\cdot\beta)}{\beta^2}\beta\in L_{2d}$, for all $x\in
L_{2d}$ and $\beta$ is primitive, $p^{e}$ divides $m$ but not
$\alpha$. Moreover $p^e$ divides $d$, because
$\beta^2=-2d\alpha^2+m^2j^2$. Thus $\beta$ satisfies (a) in
(\ref{eqn:refl}).

Notice that a reflection $r_\beta\in \OO(L_{2d})$ with $\beta^2<0$
is such that $[r_\beta]\in\mathrm{Ref}(G_{2d})$ if and only
if (A) $\frac{2d}{\beta^2}\alpha^2\not\equiv 0\pmod{d}$ and (B)
$1+\frac{2d}{\beta^2}\alpha^2\not\equiv 0\pmod{d}$. Indeed, by the
very definitions of the discriminant group $A_{L_{2d}}$ and of the
natural map $\OO(L_{2d})\to\OO(A_{L_{2d}})$, $r_\beta$ is trivial in
$G_{2d}$ if and only if either $r_\beta(k)=(2dm_1+1)k+l_1$ or
$r_\beta(k)=(2dm_2-1)k+l_2$, for $m_1,m_2\in \mathbb{Z}$ and
$l_1,l_2\in U^{\oplus 2}\oplus E_8(-1)^{\oplus 2}\hookrightarrow
L_{2d}$. This gives exactly (A) and (B).

Since $\beta^2/2$ does not divide $\alpha^2$, it follows that (A)
and (B) hold true if and only if (b) in (\ref{eqn:refl}) holds.
Hence $V(\kh_L)=\langle\beta\rangle\in\Re$ and this proves the
``if'' direction. The converse is easy and left to the
reader.\end{proof}

We discuss now two examples which show that the fixed loci of the
non-trivial elements of $G_{2d}$ may or may not contain a Heegner
divisor.

\begin{ex}\label{ex1} Suppose $2d=4p$,
where $p$ is an odd prime such that $p$ is not a square modulo 4 and
$2$ is not a square modulo $p$. By Lemma \ref{lem:act},
$G_{2d}=\langle g\rangle\iso\ZZ/2\ZZ$. Let
$\beta_1:=\alpha_1k+2j_1=\left(1,2,p-1,0,\ldots,0\right)\in
L_{2d}=\langle k\rangle\oplus U^{\oplus 2}\oplus E_8(-1)^{\oplus
2}$. As $\langle\beta_1\rangle\in\Re$, by Lemma
\ref{thm:Heegnerfixed}, $\kh_{\beta_1^\perp}\subseteq
\mathrm{Fix}(g)$ and $\beta_1^2=-4$. Moreover $\mathrm{Fix}(g)$ does
not contain any other Heegner divisor. Indeed, suppose that there
exists a primitive $\beta_2\in L_{2d}$ such that
$\kh_{\beta_2^\perp}\subseteq\mathrm{Fix}(g)$. Using Lemma
\ref{thm:Heegnerfixed}, we have $\langle\beta_2\rangle\in\Re$ and
$\beta_2^2=-4$.

Let $\beta_2=\alpha_2k+2j_2$, where $j_2\in U^{\oplus 2}\oplus
E_8(-1)^{\oplus 2}$ ($j_2$ is not necessarily primitive). Since
$-4=\beta_2^2=-2d\alpha_2^2+4j_2^2$, $\alpha_2^2\equiv 1\pmod{2}$
and there exists an integer $w$ such that $\alpha_2=1+2w$. Observe
that $\frac{1}{2}j_2\not\in U^{\oplus 2}\oplus E_8(-1)^{\oplus 2}$.
Indeed, if $\frac{1}{2}j_2\in U^{\oplus 2}\oplus E_8(-1)^{\oplus
2}$, then there would exist an integer $\alpha$ such that
$p\alpha^2\equiv 1\pmod{4}$. This is impossible since $p$ is not a
square modulo $4$.

Let $L_i$ be the minimal primitive sublattice of $\Lambda$ containing the lattice
generated by $h$ and $\beta_i$. An integral basis of $L_i$ is given
by $e_i:=h$ and $f_i:=\frac{2}{\beta^2}(\alpha_ih-\beta_i)$. We can
now define an isometry $\psi:L_2\rightarrow L_1$ such that
$e_2\mapsto e_1$ and $f_2\mapsto -we_1+f_1$. Due to \cite[Thm.\
1.14.4]{Ni}, $\psi$ extends to an isometry $\psi'\in\OO(\Lambda)$.
Hence $\kh_{\beta_1^\perp}=\kh_{\beta_2^\perp}$.

As a very special case, we have Example 1.3 in \cite{Mu1} where
Mukai's duality is studied for K3 surfaces with Picard number 1 and
with a polarization of degree 12. It is quite easy to see that, if
$\kh_L$ is the unique Heegner divisor in the fixed locus of the
non-trivial $g\in G_{2d}$ and $[\sigma]\in\kh_L$ is such that
$\rho([\sigma])=2$, then
\[
\NS(X_\sigma)\cong\left(\mathbb{Z}^2,\left(\begin{array}{cc} 2 & 4\\
4 & 2\end{array}\right)\right).
\]
In particular $\NS(X_\sigma)$ is isometric to the N{\'e}ron-Severi
group of a generic K3 surface which is a complete intersection of
bi-degree $(1,1)$ and $(2,2)$ in $\mathbb{P}^2\times\mathbb{P}^2$.
These surfaces were studied by Wehler in \cite{W}.\end{ex}

\begin{ex}\label{ex:noref} There are
$2d$-polarizations such that none of the non-trivial elements of
$G_{2d}$ fixes a Heegner divisor. An example is when $d=15$.
Indeed in this case $G_{2d}=\{\mathrm{id},g\}$ (see Lemma
\ref{lem:act}) and if there exists a Heegner divisor
$\kh_L\subseteq\mathrm{Fix}(g)$ with
$V(\kh_L)=\langle\beta\rangle$, then, by Lemma
\ref{thm:Heegnerfixed}, $\beta^2$ is either $-6$ or $-10$. If
$\beta^2=-6$, then
$-15\alpha^2+3^2s=-d\alpha^2+m^2\frac{j^2}{2}=\frac{\beta^2}{2}=-3$
for some integer $s$. In particular, $5$ should be a square modulo
$3$ but this is not the case. Similarly, $\beta^2$ cannot be equal
to $-10$. Hence $\mathrm{Fix}(g)$ does not contain a
divisor.\end{ex}

\medskip

{\small\noindent{\bf Acknowledgements.} During the preparation of
this paper the author was partially supported by the MIUR of the
Italian Government in the framework of the National Research Project
``Geometry on Algebraic Varieties'' (Cofin 2004). It is a pleasure
to thank Bert van Geemen for useful discussions and many suggestions
about the results in this paper. I wish to thank Gilberto Bini for
comments on earlier versions of the paper.}



\begin{thebibliography}{99}

\bibitem{Ge} \emph{G{\'e}ometrie des surfaces K3: modules et p{\'e}riodes}, Asterisque {\bf 126}
(1985).

\bibitem{BPV} W. Barth, K. Hulek, C. Peters, A. Van de Ven, \emph{Compact
complex surfaces}, Springer-Verlag, Berlin (2004).

\bibitem{Bo} C. Borcea, \emph{K3 surfaces and complex multiplication},
Rev. Roumaine Math. Pures Appl. {\bf 31} (1986), 499--505.

\bibitem{BM} T. Bridgeland, A. Maciocia, \emph{Complex surfaces with equivalent derived categories},
Math. Z. {\bf 236} (2001), 677--697.


\bibitem{Do} I.V. Dolgachev, \emph{Mirror symmetry for lattice polarized
K3 surfaces}, J. Math. Sci. {\bf 81} (1996), 2599--2630.


\bibitem{H2}  B. Hassett, \emph{Special cubic fourfolds}, Compositio Math. {\bf 120} (2000),
1--23.


\bibitem{HLOY1} S. Hosono, B.H. Lian, K. Oguiso, S.-T. Yau, {\it
Fourier-Mukai number of a K3 surface}, CRM Proc. and Lect. Notes
{\bf 38} (2004), 177--192.

\bibitem{Mo} D.R. Morrison, \emph{On K3 surfaces with large Picard
number}, Invent. Math. {\bf 75} (1984), 105--121.

\bibitem{Mu} S. Mukai, \emph{Duality between $D(X)$ and $D(\widehat X)$ with its applications to Picard sheaves},
Nagoya Math. J. {\bf 81} (1981), 153--175.

\bibitem{Mu2} S. Mukai, \emph{On the moduli space of bundles on K3
surfaces, I}, In: Vector Bundles on Algebraic Varieties, Oxford
University Press, Bombay and London (1987), 341--413.

\bibitem{Mu1} S. Mukai, \emph{Duality of polarized K3 surfaces}, New
trends in algebraic geometry (Warwick, 1996), London Math. Soc.
Lecture Note Ser. {\bf 264}, Cambridge Univ. Press, Cambridge
(1999), 311--326.

\bibitem{Ni} V.V. Nikulin, \emph{Integral symmetric bilinear forms and
some of their applications}, Math. USSR Izvestija {\bf 14} (1980),
103--167.

\bibitem{Ni1} V.V. Nikulin, \emph{Correspondences of a K3 surfaces with itself via moduli of sheaves. I}, math.AG/0609233.

\bibitem{Og} K. Oguiso, \emph{K3 surfaces via almost-primes}, Math.
Research Letters {\bf 9} (2002), 47--63.

\bibitem{Or} D. Orlov, \emph{Equivalences of derived categories and K3
surfaces}, J. Math. Sci. {\bf 84} (1997), 1361--1381.

\bibitem{Sa} B. Saint-Donat, \emph{Projective models of K-3 surfaces},
Amer. J. Math. {\bf 96} (1974), 602--639.

\bibitem{Sc} F. Scattone, \emph{On the compactification of moduli spaces for algebraic K3 surfaces},
Mem. Amer. Math. Soc. {\bf 70} (1987), no. 374.

\bibitem{St}  P. Stellari, \emph{Some remarks about the FM-partners
of K3 surfaces with Picard numbers 1 and 2}, Geom. Dedicata {\bf
108} (2004), 1--13.

\bibitem{Sr} H. Sterk, \emph{Finiteness results for algebraic K3
surfaces}, Math. Z. {\bf 189} (1985), 507--513.

\bibitem{W} J. Wehler, \emph{K3 surfaces with Picard number 2}, Arch. Math. {\bf 50} (1988),
73--82.


\end{thebibliography}
\end{document}